\documentclass[oneside]{amsart}
\usepackage{amsmath,amssymb,color,amsthm,graphicx,cite}
\usepackage[T1]{fontenc}
\usepackage[utf8]{inputenc}
\usepackage{color,bm}
\graphicspath{{../../img_pdf/}{../img_pdf/}{./img_pdf/}}
\usepackage{enumerate} 
\usepackage[backref]{hyperref}
\RequirePackage{luatex85}
\usepackage[all]{xy}
\newtheorem{theorem}{Theorem}[section]
\newtheorem{proposition}[theorem]{Proposition}
\newtheorem{lemma}[theorem]{Lemma}
\newtheorem{corollary}[theorem]{Corollary}
\newtheorem{claim}{Claim}

\theoremstyle{definition}
\newtheorem{definition}{Definition}
\newtheorem{main}{Theorem}

\newtheorem{main_cor}[main]{Corollary}

\def\K{\mathbb{K} } 

\def\F{\mathcal{F} }

\def\R{\mathbb{R} } 
 
\def\Z{\mathbb{Z} } 
\def\nbd{neighborhood } 
 
\def\R{\mathbb{R} }

\def\-{\ominus} 
\def\+{\oplus} 
\def\0{\circ}

\author{Tomoo Yokoyama}
\date{\today}
\address{Department of Mathematics, Faculty of Science, Saitama University, Shimo-Okubo 255, Sakura-ku, Saitama-shi, 338-8570 Japan\\}
\email{tyokoyama@rimath.saitama-u.ac.jp}
\thanks{The author was partially supported by JSPS Grant Number 24K06733}

\makeatletter
\@namedef{subjclassname@2020}{%
  \textup{2020} Mathematics Subject Classification}
\makeatother

\usepackage{color,bm}

\title[Coarse non-wandering sets \& their filtration]{Coarse non-wandering sets and their filtration}


\begin{document}

\begin{abstract}
This paper investigates recurrence properties of dynamical systems under the restriction that control is available only through inputs and outputs. We introduce the concept of ``coarse non-wandering'', a generalization of the classical non-wandering concept, and construct an associated filtration based on levels that quantify the closeness of recurrence behavior under input/output-only control. The forward direction of this filtration describes how the level of control relates to recurrence properties, whereas the backward direction captures the robustness of such behaviors and, in particular, guarantees controllability through control applied only at the observation points when the observational noise is sufficiently small. Furthermore, we demonstrate that the existence of a wandering domain is equivalent to the presence of an orbit reachable within finite error but unable to return within any slightly enlarged error bound.
\end{abstract}

\maketitle

\section{Introduction} 

\subsection{Topological background}

Topological methods have become increasingly important in the analysis of high-dimensional data, particularly through the development of filtrations in applied and computational topology. 
In particular, persistent homology provides a robust framework to extract qualitative features from noisy data by tracking the evolution of topological structures across multiple scales. 
Inspired by this perspective, this paper introduces a new topological descriptor for dynamical systems, which we term a ``non-wandering filtration'', serving as an analogue of the persistence diagram and 
opening the possibility of applying persistent homology to the analysis of dynamical systems.
Such a non-wandering filtration generalizes the classical concept of non-wandering sets by introducing a hierarchy of approximations that capture the stability and recurrence behaviors of dynamical systems. 
Similar to the persistence diagram, this construction provides a hierarchical structure that reflects how robust recurrence behaviors are under perturbations and observational noise. 
In addition to its theoretical interest, the approach can also serve as a practical tool for identifying qualitative features in dynamical data.
As an illustration, one potential application of the new framework is to numerical simulations, where the filtration may reveal robust dynamical structures from input/output-based observations.

\subsection{Background from dynamical systems}

Birkhoff introduced and studied the concepts of non-wandering points and recurrent points and investigated the asymptotic behavior of orbits \cite{birkhoff1927dynamical}.
The concept of $\varepsilon$-pseudo-orbits was introduced by Conley~\cite{conley1978isolated,conley1988gradient} and Bowen~\cite{bowen1975omega,Bowen1975,bowen1978axiom} for continuous and discrete dynamical systems, and has become a fundamental tool in the study of dynamical systems, particularly through numerical simulations.
Building on the concept of $\varepsilon$-pseudo-orbits, Conley also defined a weaker form of recurrence, known as the chain recurrence \cite{conley1988gradient}.

In \cite{yokoyama2025coarse}, the concept of $\varepsilon$-chain recurrence was recently introduced to analyze recurrent behaviors and the persistence of attractors, leading to filtrations associated with dynamical systems that quantify the strength of recurrence properties.
Furthermore, the underlying idea of such filtrations has been extended to systems with gradient-like dynamics or attractors, through a concept analogous to $\varepsilon$-pseudo-orbits \cite{imoto2025filtrations}.
These filtrations measure the robustness of attractors in gradient-like dynamics, giving lower bounds on the energy required to perturb the system into a different attractor as well as on the magnitude of energy input needed at each step, and have even been applied to numerical simulations of tropical cyclone dynamics.


\subsection{Statements of main results}

In this paper, we introduce filtrations that coarsen the concept of non-wandering behavior and provide lower bounds for simplified control schemes, which avoid computationally intensive costs such as control or perturbation at every step, by considering the non-wandering property as a property returning approximately to a neighborhood of a point. 
We now state our main results.
In fact, by introducing the concept of coarse non-wandering, we obtain the following filtrations for continuous mappings on metric spaces.

\begin{main}\label{main:01}
For any continuous mapping $f \colon X \to X$ on a metric space $X$, the family $(\Omega_{\varepsilon}(f))_{\varepsilon \in \R}$ {\rm(}see Definitions~\ref{def:cr} and \ref{def:cr02} below{\rm)} is a filtration with $\Omega(f) = \Omega_{0}(f)$ and $X = \bigcup_{\varepsilon \in \R}\Omega_{\varepsilon}(f)$. 
\end{main}

A more general statement is described in Theorem~\ref{th:main01} for any mappings on metric spaces. 
%
Furthermore, we obtain the following characterization, which leads to the application that the existence of a wandering domain is guaranteed by finding a point that can be reached within some finite error, but cannot be returned to within the error plus a small amount.

\begin{main}\label{main:02}
The following statements are equivalent for a continuous mapping $f$ on a metric space $X$: 
\\
{\rm(1)} $\Omega(f) = X$. 
\\
{\rm(2)} $\Omega_{-\varepsilon}(f) = X$ for any number $\varepsilon \geq 0$. 
\\
{\rm(3)} $z \sim^1_{\varepsilon+} x$ (see Defnition~\ref{def:e_link+} below) for any point $x \in X$, any number $\varepsilon \geq 0$, and any point $z \in [x]^1_\varepsilon$ (see Defnition~\ref{def:e_link_set} below).
\end{main}

The previous theorem implies a characterization of the existence of wandering domains as follows. 
 
 \begin{main_cor}
The following statements are equivalent for a continuous mapping $f$ on a metric space $X$: 
\\
{\rm(1)} The mapping $f$ has a wandering domain. 
\\
{\rm(2)} There are a point $x \in X$, a number $\varepsilon \geq 0$, and a point $z \in [x]^1_\varepsilon$ such that $z \notin \bigcap_{\varepsilon'>\varepsilon} [x]^1_{\varepsilon}$.
\\
{\rm(3)} There are a point $x \in X$, a number $\varepsilon \geq 0$, a point $z \in [x]^1_\varepsilon$, and a number $\varepsilon'>\varepsilon$ such that there is no $\varepsilon'$-link from $x$ to $z$ (see Defnition~\ref{def:e_link} below). 
 \end{main_cor}

In addition, one can also derive a similar result in the case of flows (see Theorem~\ref{th:02_flow} and Corollary~\ref{cor:wandering_domain_ch} below).
Moreover, we observe the following singular limit behavior of $\Omega_{-\varepsilon}$ for mappings that exist (see Theorem~\ref{th:counter_ex04}). 

The present paper consists of eight sections.
In the next section, we review some key concepts from combinatorics and dynamical systems. 
We discuss mappings in \S~3--4 and semiflows in \S~5--6. 
In fact, the $\varepsilon$-non-wandering property for mappings is introduced in \S~3, and we introduce ``non-wandering property with negative errors'' in \S~4. 
Similarly, in the context of the semiflows, we also introduce the $\varepsilon$-non-wandering property in \S~5 and ``non-wandering property with negative errors'' in \S~6. 
In \S~7, we define the non-wandering diagram to analyze dynamical systems, which is analogous to a bifurcation diagram. 
In the final section, we discuss a condition of ``negative coarse non-wandering'' property by constructing a continuous mapping.

\section{Preliminaries}

To analyze the various phenomena, we define filtration and cost function as follows.

\begin{definition}\label{def:filtration}
Let $X$ be a set and $\F = \{ F_i \mid i \in I\} \subset 2^X$ a family indexed by a totally ordered set $I$, where $2^X$ is the power set of $X$. The family $\F$ is a {\bf filtration} if $X = \bigcup_{i \in I} F_i$ and $F_{i_1} \subseteq F_{i_2}$ for any pair $i_1 \leq i_2 \in I$.
\end{definition}


\begin{definition}\label{def:metric-like}
A cost function $c \colon X^2 \to [0,\infty]$ is {\bf non-degenerate} if $c^{-1}(0) = \{ (x,x) \mid x \in X\}$. 
\end{definition}

Notice that a metric is a non-degenerate cost function. 
This means that the conventional framework of dynamical systems on metric spaces is generalized to a broader framework of dynamical systems on sets equipped with a cost function.


\subsection{Basic concepts in mappings}
Let $f \colon X \to X$ be a mapping on a topological space $X$, which need not be continuous. 
For any point $x \in X$, we denote $O^+(x)$ the positive orbit (i.e. $O^+(x) := \{ f^n(x) \mid n \in \Z_{> 0} \}$). 
A point $x \in X$ is {\bf fixed} if $f(x) = x$, and is {\bf periodic} if there is a positive integer $n$ with $f^n(x)= x$.
%
We define the non-wandering property for mappings as follows. 

\begin{definition}
A point is {\bf wandering} if there are its neighborhood $U$ and a positive number $N$ such that $f^n(U) \cap U = \emptyset$ for any $n > N$. 
\end{definition}

\begin{definition}
A point is {\bf non-wandering} if it is not wandering (i.e. for any its neighborhood $U$, there is a number $n \in \mathbb{Z}_{>0}$ such that $f^n(U) \cap U \neq \emptyset$).
\end{definition}

Denote by $\Omega (f)$ the set of non-wandering points, called the {\bf non-wandering set}.

\subsection{Basic concepts for homeomorphisms}

Let $f \colon X \to X$ be a bijection on a topological space $X$.
For a point $x \in X$, we denote by $O(x)$ the orbit of $x$ by $f$ and $O^-(x)$ the negative orbit (i.e. $O^-(x) := \{ f^n(x) \mid n \in \Z_{< 0} \}$).
%

\subsection{Basic concepts for semiflows}

A {\bf flow} is a continuous $\R$-action on a topological space.
A {\bf semiflow} is a continuous $\R_{\geq 0}$-action on a topological space.

Let $v$ be a continuous action $v \colon \K \times X \to X$, where $\K$ is either $\R$ or $\R_{\geq 0}$. 
Notice that if $\K = \R$ (resp. $\K = \R_{\geq 0}$) then $v$ is a flow (resp. semiflow). 
For $t \in \K$, define $v^t : X \to X$ by $v^t := v(t, \cdot )$.

A point $x \in X$ is {\bf singular} if $x = v^t(x)$ for any $t \in \K$. 
For a point $x \in X$, we denote by $O(x)$ the orbit of $x$ by $v$ (i.e. $O(x) := \{ v^t(x) \mid t \in \K \} = v(\K,x)$), and the positive orbit (i.e. $O^+(x) := \{ v^t(x) \mid t > 0 \} = v(\R_{>0},x)$). 
%
%
We define the non-wandering property as follows. 

\begin{definition}
A point is {\bf wandering} if there are its neighborhood $U$ and a positive number $N$ such that $v^t(U) \cap U = \emptyset$ for any $t > N$. 
\end{definition}

\begin{definition}
A point is {\bf non-wandering} if it is not wandering (i.e. for any its neighborhood $U$ and for any positive number $N$, there is a number $t \in \mathbb{R}$ with $t > N$ such that $v^t(U) \cap U \neq \emptyset$).
\end{definition}

Denote by $\Omega (v)$ the set of non-wandering points, called the {\bf non-wandering set}. 

\section{Coarse non-wandering property for mappings} 

In this section, we introduce the $\varepsilon$-non-wandering property. 
In this section, from now on, let $f \colon X \to X$ be a mapping on a set $X$ with a cost function $c \colon X^2 \to [0,\infty]$ unless otherwise stated. 
Moreover, when $X$ is a metric space equipped with a metric $d$, we assume that its cost function $c$ coincides with the metric $d$. 

\subsection{$\varepsilon$-non-wandering property}

To define the $\varepsilon$-non-wandering property, we introduce some concepts and characterize the non-wandering property. 
First, we define $\varepsilon$-links for mappings, a binary relation, and a subset as follows. 

\begin{definition}\label{def:e_link}
For any number $\varepsilon \geq 0$, define a binary relation $\sim^1_\varepsilon$ on $X$ by $x \sim^1_\varepsilon y$ if there are a piont $z \in X$ and a natural number $n \in \Z_{>0}$ such that $c(x,z) \leq \varepsilon$ and $c(f^n(z),y) \leq \varepsilon$. 
Then the subset $(f^i(z))_{i=0}^n$ is called an {\bf $\bm{\varepsilon}$-link} from $x$ to $y$ with respect to $c$. 
\end{definition}

\begin{definition}\label{def:e_link_set}
Define a subset $[x]^1_\varepsilon$ as follows: 
\[
[x]^1_\varepsilon := \{ y \in X \mid x \sim^1_\varepsilon y \} = \{ y \in X \mid \text{There is an } \varepsilon \text{-link from }x \text{ to } y \}	
\]
\end{definition}

We have the following observations. 

\begin{lemma}\label{lem:001}
For any non-negative numbers $\varepsilon_1 < \varepsilon_2$, we have $[x]^1_{\varepsilon_1} \subseteq [x]^1_{\varepsilon_2}$. 
\end{lemma}

\begin{lemma}\label{lem:cost_orbit}
If the cost function $c$ is non-degenerate, then $[x]^1_0 = O^+(x)$ for any $x \in X$.
\end{lemma}


We define a binary relation and a subset as follows. 

\begin{definition}
For any $\varepsilon \geq 0$, we define a binary relation $\sim^1_{\varepsilon +}$ on $X$ as follows: 
\[
x \sim^1_{\varepsilon +} y \text{ if } x \sim^1_{\varepsilon'} y \text{ for any } \varepsilon' > \varepsilon 
\]
\end{definition}

\begin{definition}\label{def:e_link+}
Define a subset $[x]^1_{\varepsilon +}$ as follows: 
\[
\begin{split}
[x]^1_{\varepsilon +} := &\{ y \in X \mid x \sim^1_{\varepsilon +} y \} 
\\
= &\{ y \in X \mid \text{For any }\varepsilon' > \varepsilon, \text{ there is an } \varepsilon' \text{-link from }x \text{ to } y \}	
\end{split}
\]
\end{definition}

Then we have the following observation. 


\begin{lemma}\label{lem:02}
We have the following statements and for any $x \in X$: 
\\
{\rm(1)} For any non-negative number $\varepsilon \geq 0$, we have $[x]^1_{\varepsilon +} = \bigcap_{\varepsilon' > \varepsilon} [x]^1_{\varepsilon'}$. 
\\
{\rm(2)} For any non-negative numbers $\varepsilon_1 \leq \varepsilon_2$, we have $[x]^1_{\varepsilon_1 +} \subseteq [x]^1_{\varepsilon_2 +}$. 
\end{lemma}

We characterize the non-wandering property for mappings.

\begin{lemma}\label{lem:nw_eq}
Let $f \colon X \to X$ be a mapping on a metric space $(X,d)$. 
Then the following statements are equivalent for any point $x \in X$:
\\
{\rm(1)} The point $x$ is non-wandering. 
\\
{\rm(2)} For any $\varepsilon > 0$, there is an $\varepsilon$-link from $x$ to $x$ with respect to the cost function $d$ {\rm(i.e.} $x \sim^1_{0+} x$ with respect to $d${\rm)}. 
\end{lemma}

\begin{proof}
Notice that $x$ is non-wandering if and only if, for any $\varepsilon > 0$, there is a number $n \in \mathbb{Z}_{>0}$ such that $f^n(B_\varepsilon (x)) \cap B_\varepsilon (x) \neq \emptyset$, where $B_\varepsilon (x) := \{ y \in X \mid d(x,y) < \varepsilon \}$ is the open $\varepsilon$-ball centered at $x$. 
Since the condition $f^n(B_\varepsilon (x)) \cap B_\varepsilon (x) \neq \emptyset$ is equivalent to the condition $x \sim^1_\varepsilon x$ with respect to $d$, the point $x$ is non-wandering if and only if $x \sim^1_\varepsilon x$ for any $\varepsilon > 0$. 
\end{proof}

By the characterization in the previous lemma, we introduce the $\varepsilon$-non-wandering property for mappings as follows. 

\begin{definition}\label{def:cr}
A point $x \in X$ is {\bf $\bm{\varepsilon}$-non-wandering} if $x \sim^1_{\varepsilon +} x$. 
\end{definition}

On other words, a point $x \in X$ is $\varepsilon$-non-wandering if and only if, for any $\varepsilon' > \varepsilon$, there is a point $y \in X$ with $c(x,y) \leq \varepsilon'$ and there is a positive integer $n \in \Z_{>0}$ such that $c(f^n(y),x) \leq \varepsilon'$. 

\begin{definition}
The set of non-wandering points is called the {\bf $\bm{\varepsilon}$-non-wandering set} of $f$ and is denoted by $\bm{\Omega_\varepsilon(f)}$. 
\end{definition}

Notice that the $0$-non-wandering property with respect to the metric corresponds to the original wandering property. 
%


\subsection{Properties of the $\varepsilon$-non-wandering property}

%




We have the following observations. 

\begin{lemma}\label{lem:002}
We have the following statements for any non-negative numbers $\varepsilon_1 < \varepsilon_2 < \varepsilon_3$ and for any $x \in X$: 
\[
[x]^1_{\varepsilon_1 +} \subseteq [x]^1_{\varepsilon_2} \subseteq [x]^1_{\varepsilon_2 +} \subseteq [x]^1_{\varepsilon_3}
\]
\end{lemma}

\begin{proof}
Lemma~\ref{lem:001} implies $[x]^1_{\varepsilon_2} \subseteq [x]^1_{\varepsilon_2 +}$.
By definitions, we have $[x]^1_{\varepsilon_2 +} = \{ y \in X \mid x \sim^1_{\varepsilon_2 +} y \} \subseteq \{ y \in X \mid x \sim^1_{\varepsilon_3} y \} = [x]^1_{\varepsilon_3 +}$. 
The same argument implies $[x]^1_{\varepsilon_1 +} \subseteq [x]^1_{\varepsilon_2}$. 
\end{proof}
%

%


Notice that $\Omega(f) = \Omega_{0}(f) \subseteq \Omega_{\varepsilon}(f)$ for any $\varepsilon \geq 0$ if the cost function $c$ is a metric on $X$. 
We obtain the following observation. 

\begin{proposition}\label{lem:equality}
Let $f \colon X \to X$ be a mapping on a set $X$ with a cost function $c \colon X^2 \to [0,\infty]$. 
The family $(\Omega_{\varepsilon}(f))_{\varepsilon \in [0,\infty]}$ is a filtration. 
In fact, the following properties hold: 
\\
{\rm(1)} For any non-negative numbers $\varepsilon_1 \leq \varepsilon_2$, we have $\Omega_{\varepsilon_1 }(f) \subseteq \Omega_{\varepsilon_2}(f)$. 
\\
{\rm(2)} $X = \bigcup_{\varepsilon \in [0,\infty]} \Omega_{\varepsilon}(f)$. 
\\
{\rm(3)} If $X$ is a metric space whose cost function is the metric, then $X = \bigcup_{\varepsilon \in \R_{\geq 0}} \Omega_{\varepsilon}(f)$. 
\end{proposition}

\begin{proof}
By definition of $\Omega_{\varepsilon}(f)$, assertion (1) holds. 
For any $x \in X$, we have $x \in \Omega_{c(f(x),x)}(f)$ and so assertion (2) holds. 

Suppose that $X$ is a metric space whose cost function is the metric $d$. 
Then $d^{-1}(\infty) = \emptyset$. 
Assertion (3) follows from assertion (2). 
\end{proof}

%

%

\section{$-\varepsilon$-non-wandering property for mappings}

To extend the filtration $(\Omega_{\varepsilon}(f))_{\varepsilon \geq 0}$, we introduce an $-\varepsilon$-non-wandering property for mappings as follows. 
Let $f \colon X \to X$ be a mapping on a set $X$ with a cost function $c \colon X^2 \to [0,\infty]$. 
We define the following binary relation to formulate ``non-wandering property with negative errors''. 

\begin{definition}
For any $\varepsilon \geq 0$, we define a binary relation $\sim^1_{-\varepsilon +}$ as follows: 
\[
x \sim^1_{-\varepsilon +} y \text{ if } z \sim^1_{\varepsilon' +} y \text{ for any } \varepsilon' \in [0,\varepsilon] \text{ and for any } z \in [x]^1_{\varepsilon'}
\]
\end{definition}

We introduce the following concept. 

\begin{definition}\label{def:cr02}
For any non-negative number $\varepsilon \geq 0$, a point $x \in \Omega_0(f)$ is {\bf $\bm{-\varepsilon}$-non-wandering} if $x \sim^1_{-\varepsilon +} x$. 
Denote by $\bm{\Omega_{-\varepsilon}(f)}$ the set of $-\varepsilon$-non-wandering points. 
\end{definition}

Notice that the condition $x \sim^1_{-0 +} x$ need not imply $x \in \Omega_0(f)$ (see an example in Lemma~\ref{lem:couter_example01} in details). 
On the other hand, the equivalence holds for homeomorphisms. 

\begin{lemma}\label{lem:41}
Let $f \colon X \to X$ be a homeomorphism on a metric space $(X,d)$ whose cost function is the metric. 
For any $x \in X$ with $x \sim^1_{0 +} x$, we have $x \in \Omega_{0}(f)$. 
\end{lemma}

\begin{proof}
Fix any point $x \in X$ with $x \sim^1_{0 +} x$ and any positive number $\varepsilon > 0$. 
By the continuity of $f^{-1}$, there is a positive number $\delta > 0$ such that $f^{-1}(B_\delta(f(x))) \subseteq B_\varepsilon (x)$. 
Pur $\varepsilon_1 := \min \{ \varepsilon, \delta \} > 0$. 
Since $f(x) \in [x]^1_0$, by $x \in \Omega_{-0}(f)$, there are a point $z \in B_{\varepsilon_1} (f(z))$ and a nutural number $n \in \Z_{>0}$ such that $d(f^n(z),x) < \varepsilon_1 \leq \varepsilon$. 
From $z \in B_{\varepsilon_1} (f(z)) \subseteq B_\delta(f(x))$, we have $f^{-1}(z) \in B_\varepsilon (x)$. 
Then the sequence $(f^{-1}(z), z, f(z), \ldots , f^n(z))$ is an $\varepsilon$-link from $x$ to $x$. 
This means that $x \in \Omega(f) = \Omega_{0}(f)$. 
Therefore, we obtain $\Omega_{-0}(f) \subseteq \Omega_{0}(f)$.
%
\end{proof}

By definition, the following statement holds. 

\begin{lemma}\label{lem:41-}
For any non-negative numbers $\varepsilon_1 < \varepsilon_2$, we have $\Omega_{-\varepsilon_2}(f) \subseteq \Omega_{-\varepsilon_1}(f)$. 
\end{lemma}

%

We have the following persistence. 

\begin{lemma}\label{lem:equv_minus}
Let $f \colon X \to X$ be a continuous mapping on a metric space $X$. 
If $\Omega(f) = X$, then $\Omega_{-\varepsilon}(f) = X$ for any non-negative number $\varepsilon \geq 0$. 
\end{lemma}

\begin{proof}
Fix any point $x \in X$ and any non-negative number $\varepsilon \geq 0$. 
Choose any $\varepsilon' \in [0,\varepsilon]$, any $y \in [x]^1_{\varepsilon'}$, and any $\varepsilon_1>0$. 
Put $\varepsilon'' := \varepsilon' + \varepsilon_1>0$. 
It suffices to show  that $y \sim^1_{\varepsilon''} x$. 
By $y \in [x]^1_{\varepsilon'}$, there is an $\varepsilon'$-link $(f^i(z))_{i=0}^n$ from $x$ to $y$. 
Then $d(x,z) <  \varepsilon'$ and $d(f^n(z),y) < \varepsilon'$. 
From the continuity of $f^n$, there is a positive number $\delta>0$ such that $f^n(B_\delta(z)) \subseteq B_{\varepsilon_1}(f^n(z))$. 
Put $\varepsilon_2 := \min \{ \delta, \varepsilon_1 \}$. 
Since $z$ is non-wandering, there are a point $z' \in B_{\varepsilon_2}(z)$ and a natural number $n' > n$ such that $d(f^{n'}(z'), z) < \varepsilon_2$.  
Then $d(y,f^n(z')) \leq d(y,f^n(z)) + d(f^n(z),f^n(z')) < \varepsilon' + \varepsilon_1 = \varepsilon''$ and $d(f^{n'}(z'), x) \leq d(f^{n'}(z'), z) + d(z,x) < \varepsilon_2 + \varepsilon' \leq \varepsilon''$. 
Therefore, a sequence $(f^i(z))_{i=n}^{n'}$ is an $\varepsilon''$-link from $y$ to $x$. 
This means that $y \sim^1_{\varepsilon''} x$. 
\end{proof}

We now prove the following as one of our main results. 

\begin{proof}[Proof of Theorem~\ref{main:02}]
The definition of $\varepsilon$-non-wandering point implies that assertions (2) and (3) are equivalent. 
By Lemma~\ref{lem:equv_minus}, assertions (1) and (2) are equivalent. 
%
 \end{proof}
  
 \subsection{Existence of filtrations for coarse non-wandering property}
 
We have the following correspondence. 

\begin{lemma}\label{lem:-zero_zero}
Let $f$ be a continuous mapping on a metric space $(X,d)$.
Then $\Omega_{-0}(f) = \Omega_{0}(f) = \Omega(f)$.
\end{lemma}

\begin{proof}
By definitions, we have $\Omega_{-0}(f)  \subseteq \Omega_0(f) = \Omega(f)$. 
Therefore, we may assume that $\Omega_{0}(f) \neq \emptyset$. 
Fix any point $x \in \Omega_{0}(f)$. 
Take any point $z \in [x]^1_{0}$ and any $\varepsilon > 0$. 

\begin{claim}
 $z \sim^1_{\varepsilon} x$. 
\end{claim}

\begin{proof}
Since $d$ is the cost fuction, Lemma~\ref{lem:cost_orbit} implies that $[x]^1_{0} = O^+(x)$. 
Then there is a positive integer $n \in \Z_{>0}$ such that $f^n(x) = z$. 
By the continuity of $f$, there is a positive number $\delta \in (0,\varepsilon)$ such that $d(f^n(x'),f^n(x)) < \varepsilon$ for any $x' \in X$ with $d(x,x') < \delta$. 
From $x \in \Omega_0(f)$, there are a point $x' \in X$, a positive integer $N > n$, and a $\delta$-link $(f^i(x'))_{i=0}^N$ from $x$ to $x$. 
Then we have the following inequalities: 
\[
d(f^n(x),f^n(x')) < \varepsilon
\]
\[
d(f^N(x'),x) < \delta < \varepsilon
\]
This means that the sequence $(f^i(x'))_{i=n}^N$ is an $\varepsilon$-link from $z = f^n(x)$ to $x$. 
Therefore, we have $z \sim^1_{\varepsilon} x$. 
\end{proof}
Since $\varepsilon > 0$ can be arbitrarily small, we obtain $z \sim^1_{0+} x$. 
This means that $x \in \Omega_{-0}(f)$. 
\end{proof}

\subsubsection{Total order on a variance of $\R$}

Consider the pair $(\R, \leq_{\R})$ of the real line and the standard order. 
Add a maximal point $\infty$ to $\R$ and the resulting torally ordered set is denoted by $(-\infty, \infty]$. 
Write two distinct points $-0, +0$ which are not contained in $\R$. 
Set $(-\infty,-0] := (\infty,0) \sqcup \{ -0 \}$ and $[0,\infty] := \{ +0 \} \sqcup (0,\infty]$. 
By considering $-0$ as the maximal element in $(-\infty,-0]$ and $+0$ as the minimal element in $[0,\infty]$, the subsets $(-\infty,-0]$ and $[0,\infty]$ are totally ordered. 
By setting $-0<+0$, the disjoint union $(-\infty,-0] \sqcup [0,\infty]$ becomes a totally ordered set. 

\subsubsection{Filtration whose parameters forms $(-\infty,-0] \sqcup [0,\infty]$}

Proposition~\ref{lem:equality} and the previous lemma imply the following inclusions. 

\begin{theorem}\label{th:main01}
Let $f \colon X \to X$ be a mapping on a set $X$ with a cost function $c \colon X^2 \to [0,\infty]$.  
The family $(\Omega_{\varepsilon}(f))_{\varepsilon \in (-\infty,-0] \sqcup [0,\infty]}$ of subsets of $X$ is a filtration of $X$. 
\end{theorem}

Lemma~\ref{lem:-zero_zero} and the previous theorem imply Theorem~\ref{main:01}. 
We call the families in Theorem~\ref{th:main01} and Theorem~\ref{main:01} the {\bf non-wandering filtrations} of the mappings.

\section{Coarse non-wandering property for (semi)flows}


In the following sections, we demonstrate that some results analogous to those obtained for the mapping also hold for semiflows. 
In this section, we introduce $\varepsilon$-non-wandering property for semiflows. 

Let $v$ be a continuous action $v \colon \K \times X \to X$ on a topological space $X$ with a cost function $c \colon X^2 \to [0,\infty]$, where $\K$ is either $\R$ or $\R_{\geq 0}$. 
Notice that if $\K = \R$ (resp. $\K = \R_{\geq 0}$) then $v$ is a flow (resp. semiflow). 
Recall that, when $X$ is a metric space equipped with a metric $d$, we assume that the cost function $c$ coincides with the metric $d$.

\subsection{Characterization of non-wandering property for (semi)flows}

First, we define $(\varepsilon,T)$-links for semiflows and a binary relation as follows. 

\begin{definition}\label{def:e_link_flow}
For any number $\varepsilon \geq 0$ and any number $T>0$, define a binary relation $\sim^1_{(\varepsilon,T)}$ on $X$ by $x \sim^1_{(\varepsilon,T)} y$ if there are a piont $z \in X$ and a number $r \in \R_{\geq T}$ such that $c(x,z) \leq \varepsilon$ and $c(v^r(z),y) \leq \varepsilon$. 
Then the subset $(v^t(z))_{t \in [0,r]}$ is called an {\bf $\bm{(\varepsilon,T)}$-link} from $x$ to $y$ with respect to $c$. 
\end{definition}

\begin{definition}
Define $[x]^1_{(\varepsilon,T)}$ as follows: 
\[
[x]^1_{(\varepsilon,T)} := \{ y \in X \mid x \sim^1_{(\varepsilon,T)} y \} = \{ y \in X \mid \text{ there is an } (\varepsilon,T) \text{-link from }x \text{ to } y \}
\]  
\end{definition}

\begin{definition}
Define the following binary relation $\sim^1_\varepsilon$ for any $\varepsilon \geq 0$:  
\[
x \sim^1_{\varepsilon} y \text{ if } x \sim^1_{(\varepsilon,T)} y \text{ for any } T > 0
\]
\end{definition}

\begin{definition}
Define $[x]^1_{\varepsilon}$ as follows: 
\[
\begin{split}
[x]^1_{\varepsilon} := & \{ y \in X \mid x \sim^1_{\varepsilon} y \} 
\\
= & \{ y \in X \mid x \sim^1_{(\varepsilon,T)} y \text{ for any } T > 0 \}
\\
= & \{ y \in X \mid \text{For any } T > 0, \text{ there is an } (\varepsilon,T) \text{-link from }x \text{ to } y \}
\end{split}
\] 
\end{definition}

Notice that $[x]^1_{\varepsilon} = \bigcap_{T > 0} [x]^1_{(\varepsilon,T)}$. 
We have the following observation. 

\begin{lemma}
For any non-negative numbers $\varepsilon_1 < \varepsilon_2$, we have $[x]^1_{\varepsilon_1} \subseteq [x]^1_{\varepsilon_2}$. 
\end{lemma}

We define a binary relation and a subset as follows. 

\begin{definition}
For any $\varepsilon \geq 0$, we define a binary relation $\sim^1_{\varepsilon +}$ on $X$ as follows: 
\[
x \sim^1_{\varepsilon +} y \text{ if } x \sim^1_{\varepsilon'} y \text{ for any } \varepsilon' > \varepsilon 
\]
\end{definition}

\begin{definition}\label{def:e_link_flow+}
Define a subset $[x]^1_{\varepsilon +}$ as follows: 
\[
\begin{split}
[x]^1_{\varepsilon +} := &\{ y \in X \mid x \sim^1_{\varepsilon +} y \} 
\\
= & \{ y \in X \mid x \sim^1_{\varepsilon'} y \text{ for any } \varepsilon' > \varepsilon \}
\\
= & \{ y \in X \mid x \sim^1_{(\varepsilon',T)} y \text{ for any } \varepsilon' > \varepsilon \text{ and any } T>0 \}
\\
= &\{ y \in X \mid \text{For any }\varepsilon' > \varepsilon \text{ and any } T > 0, \text{ there is an } (\varepsilon',T) \text{-link from }x \text{ to } y \}	
\end{split}
\]
\end{definition}

%
%

We have the following observations. 

\begin{lemma}\label{lem:02++}
We have the following statements: 
\\
{\rm(1)} For any non-negative number $\varepsilon$, we have $[x]^1_{\varepsilon +} = \bigcap_{\varepsilon' > \varepsilon} [x]^1_{\varepsilon'} = \bigcap_{\varepsilon' > \varepsilon} \bigcap_{T > 0} [x]^1_{(\varepsilon',T)}$. 
\\
{\rm(2)} For any non-negative numbers $\varepsilon_1 < \varepsilon_2$, we have $[x]^1_{\varepsilon_1 +} \subseteq [x]^1_{\varepsilon_2 +}$. 
\end{lemma}

\subsubsection{Characterization of non-wandering property for (semi)flows}

We characterize the non-wandering property for the continuous actions.

\begin{lemma}\label{lem:nw_eq_flow}
Let $v$ be a continuous action $v \colon \K \times X \to X$ on a metric space $(X,d)$, where $\K$ is either $\R$ or $\R_{\geq 0}$. 
The following statements are equivalent for any point $x \in X$:
\\
{\rm(1)} The point $x$ is non-wandering. 
\\
{\rm(2)} For any positive numbers $\varepsilon, T >0$, there is an $(\varepsilon,T)$-link from $x$ to $x$ with respect to the cost function $d$ {\rm(i.e.} $x \sim^1_{(\varepsilon,T)} x$ with respect to $d${\rm)}. 
\\
{\rm(3)} $x \sim^1_{0+} x$ with respect to the cost function $d$. 
\end{lemma}

\begin{proof}
We have the following equivalence: 
\[
\begin{split}
& x \text{ is non-wandering}
\\
\Longleftrightarrow \,\, & \forall \varepsilon > 0, \forall T > 0, \exists t \in \mathbb{R}_{> T} \text{ s.t. } v^t(B_\varepsilon (x)) \cap B_\varepsilon (x) \neq \emptyset
\\
\Longleftrightarrow \,\, & \forall \varepsilon> 0, \forall T > 0, \exists t \in \mathbb{R}_{\geq T} \text{ s.t. } v^t(B_\varepsilon (x)) \cap B_\varepsilon (x) \neq \emptyset
\\
\Longleftrightarrow \,\, &  \forall \varepsilon> 0, \forall T > 0, \exists t \in \mathbb{R}_{\geq T}, \exists z \in B_\varepsilon (x) \text{ s.t. } v^t(z) \in B_\varepsilon (x)
\\ 
\Longleftrightarrow \,\, &  \forall \varepsilon> 0, \forall T > 0, \exists t \in \mathbb{R}_{\geq T}, \exists z \in X \text{ s.t. } d(x,z) \leq \varepsilon \text{ and }  d(v^t(z),x) \leq \varepsilon
\\ 
\Longleftrightarrow \,\, &  \forall \varepsilon> 0, \forall T > 0 : 
x \sim^1_{(\varepsilon,T)} x 
\text{ with respect to } d 
\\ 
\Longleftrightarrow \,\, &  \forall \varepsilon> 0 : x \sim^1_{\varepsilon} x \text{ with respect to } d 
\\ 
\Longleftrightarrow \,\, &   x \sim^1_{0+} x \text{ with respect to } d 
\end{split}
\]
This implies the assertion. 
\end{proof}

\subsection{$\varepsilon$-non-wandering property for (semi)flows}

By the previous lemma, we introduce the following variant of the non-wandering property and relative concepts. 

\begin{definition}
For any $\varepsilon \geq 0$, a point $x \in X$ is {\bf $\varepsilon$-non-wandering} if $x \sim^1_{\varepsilon +} x$ {\rm(i.e.} $x \in [x]^1_{\varepsilon +}${\rm)}. 
Denote by $\bm{\Omega_{\varepsilon}(v)}$ the set of $\varepsilon$-non-wandering points. 
\end{definition}


%
Moreover, we have the following observation. 

\begin{lemma}
The family $(\Omega_{\varepsilon}(v))_{\varepsilon \geq 0}$ is a filtration. 
In particular, for any non-negative numbers $\varepsilon_1 < \varepsilon_2$, we have $\Omega_{\varepsilon_1}(v) \subseteq \Omega_{\varepsilon_2}(v)$. 
\end{lemma}

Notice that the following equality in Proposition~\ref{lem:equality}(2) 
\[
X = \bigcup_{\varepsilon \geq 0} \Omega_{\varepsilon}(v)
\]
need not hold for any (semi)flow $v$ in general.  
In fact, the (semi)flow $v \colon \K \times \R \to \R$ by $v(t,x) = x+t$ on $\R$ has no $\varepsilon$-non-wandering points for any $\varepsilon \in \R_{\geq 0}$. 
%
%

\subsection{Variants of non-wandering properties for (semi)flows}

We define a binary relation and a subset as follows. 

\begin{definition}
For any $\varepsilon \geq 0$ and any $T>0$, we define a binary relation $\sim^1_{(\varepsilon +,T)}$ on $X$ as follows: 
\[
x \sim^1_{(\varepsilon +,T)} y \text{ if } x \sim^1_{(\varepsilon',T)} y \text{ for any } \varepsilon' > \varepsilon 
\]
\end{definition}

\begin{definition}\label{def:e_link_flow_+}
Define a subset $[x]^1_{\varepsilon +}$ as follows: 
\[
\begin{split}
[x]^1_{(\varepsilon +,T)} := &\{ y \in X \mid x \sim^1_{(\varepsilon +,T)} y \} 
\\
= & \{ y \in X \mid x \sim^1_{(\varepsilon',T)} y \text{ for any } \varepsilon' > \varepsilon \}
\\
= &\{ y \in X \mid \text{For any }\varepsilon' > \varepsilon, \text{ there is an } (\varepsilon',T) \text{-link from }x \text{ to } y \}	
\end{split}
\]
\end{definition}

We have the following observation. 

\begin{lemma}\label{lem:corr_T}
For any non-negative number $\varepsilon$, we have $[x]^1_{\varepsilon +} = \bigcap_{T > 0} [x]^1_{(\varepsilon+,T)}$. 
\end{lemma}

\section{$-\varepsilon$-non-wandering property for (semi)flows}

We have the following variant of the $-\varepsilon$-non-wandering property. 

Let $v$ be a continuous action $v \colon \K \times X \to X$ on a metric space $X$, where $\K$ is either $\R$ or $\R_{\geq 0}$. 
%
We introduce the following concepts, which are ``non-wandering properties with negative errors'' for  continuous actions. 

\begin{definition}
For any $\varepsilon \geq 0$, we define a binary relation $\sim^1_{(-\varepsilon +,T)}$ as follows: 
\[
x \sim^1_{(-\varepsilon +,T)} y \text{ if } z \sim^1_{(\varepsilon' +,T)} y \text{ for any } \varepsilon' \in [0,\varepsilon] \text{ and for any } z \in [x]^1_{(\varepsilon',T)}
\]
\end{definition}

\begin{definition}
For any $\varepsilon \geq 0$, we define a binary relation $\sim^1_{-\varepsilon +}$ as follows: 
\[
x \sim^1_{-\varepsilon +} y \text{ if } x \sim^1_{(-\varepsilon +,T)} y \text{ for any } T>0
\]
\end{definition}

\begin{definition}\label{def:e_link_T}
Define a subset $[x]^1_{(-\varepsilon +,T)}$ as follows: 
\[
\begin{split}
[x]^1_{(-\varepsilon +,T)} := &\{ y \in X \mid x \sim^1_{(-\varepsilon +,T)} y \} 
\\
= & \{ y \in X \mid z \sim^1_{(\varepsilon' +,T)} y \text{ for any } \varepsilon' \in [0,\varepsilon] \text{ and for any } z \in [x]^1_{\varepsilon'} \}
\\
= & \{ y \in X \mid y \in [z]^1_{(\varepsilon'+,T)} \text{ for any } \varepsilon' \in [0,\varepsilon] \text{ and for any } z \in [x]^1_{\varepsilon'} \}
\\
= & \left\{ y \in X \middle| y \in \bigcap_{z \in [x]^1_{\varepsilon'}} [z]^1_{(\varepsilon'+,T)} \text{ for any } \varepsilon' \in [0,\varepsilon] \right\}
\\
= & \bigcap_{\varepsilon' \in [0,\varepsilon]} \bigcap_{z \in [x]^1_{\varepsilon'}} [z]^1_{(\varepsilon'+,T)}
\\
= & \bigcap_{\varepsilon' \in [0,\varepsilon]} \bigcap_{\varepsilon''>\varepsilon'} \bigcap_{z \in [x]^1_{\varepsilon'}} [z]^1_{(\varepsilon'',T)}
\end{split}
\]
\end{definition}

\begin{definition}\label{def:e_link+-}
Define a subset $[x]^1_{-\varepsilon +}$ as follows: 
\[
\begin{split}
[x]^1_{-\varepsilon +} := \{ y \in X \mid x \sim^1_{-\varepsilon +} y \} = & \bigcap_{T>0} \bigcap_{\varepsilon' \in [0,\varepsilon]} \bigcap_{z \in [x]^1_{\varepsilon'}} [z]^1_{(\varepsilon'+,T)}
\\
= & \bigcap_{T>0} \bigcap_{\varepsilon' \in [0,\varepsilon]} \bigcap_{\varepsilon''>\varepsilon'} \bigcap_{z \in [x]^1_{\varepsilon'}} [z]^1_{(\varepsilon'',T)}
\end{split}
\]
\end{definition}

We introduce the following concepts as mappings, each of which is ``non-wandering property with negative errors''.

\begin{definition}
For any number $\varepsilon \geq 0$, a point $x \in \Omega_{0}(v)$ is {\bf $-\varepsilon$-non-wandering} if $x \sim^1_{-\varepsilon +} x$.
Denote by $\bm{\Omega_{-\varepsilon}(v)}$ the set of $-\varepsilon$-non-wandering points. 
\end{definition}

Taking a suspension of an example in Lemma~\ref{lem:couter_example01}, notice that the condition $x \sim^1_{-0 +} x$ need not imply $x \in \Omega_0(v)$. 
The continuity and compactness imply the following statements. 

\begin{lemma}
Let $v$ be a flow on a compact metric space. 
Then $\Omega_{-0}(v) = \Omega_{0}(v)$.
\end{lemma}

\begin{proof}
Fix any point $x \in \Omega(v)$. 
Then we have the following observation: 
\[
[x]^1_{(0,T)} = \{ y \in X \mid \text{There is a } (0,T)\text{-link from }x \text{ to } y \} = v([T,\infty),x)
\]
We have the following equivalence relations: 
\[
\begin{split}
x \in \Omega_{-0}(v) & \Longleftrightarrow x \sim^1_{-0 +} x
\\
& \Longleftrightarrow  x \sim^1_{(-0 +,T)} x \text{ for any } T>0
\\
& \Longleftrightarrow  z \sim^1_{(0 +,T)} x \text{ for any } T>0 \text{ and for any } z \in [x]^1_{(0,T)} = v([T,\infty),x)
\\
& \Longleftrightarrow  v^s(x) \sim^1_{(0 +,T)} x \text{ for any } T>0 \text{ and for any } s \geq T
\\
& \Longleftrightarrow  v^s(x) \sim^1_{(\varepsilon,T)} x \text{ for any } \varepsilon,T>0 \text{ and for any } s \geq T
\end{split}
\]
Fix any $\varepsilon,T> 0$ and any $s \geq T$. 
Put $z = v^s(x)$. 
Then it suffices to show that there is an $(\varepsilon,T)$-link from $z$ to $x$.
Since $v$ is a flow, there is a closed transverse interval $T \subset B_\varepsilon(x)$ whose interior contains $x$ such that $v(s,T) \subset B_\varepsilon(z)$. 
Then there is a positive number $\delta > 0$ such that $v([-\delta, \delta], T)$ is a \nbd of $x$ contained in $B_\varepsilon(x)$. 
By $x \in \Omega(v)$, there is a point $x' \in T$ and a number $t' > 2s$ such that $v^{t'}(x') \in B_\varepsilon(x)$. 
Put $z' := v^s(x')$ and $r:= t' -s > s \geq T$. 
From $z' = v^s(x') \in v(s,T) \subset B_\varepsilon(z)$, the subset $(v^t(z'))_{t\in [0,t' -s]}$ is an $(\varepsilon,T)$-link from $z$ to $x$. 
\end{proof}

In fact, we have the following persistence.  

\begin{lemma}
Let $v$ be a flow on a metric space $X$.
If $\Omega(v) = X$, then $\Omega_{-\varepsilon}(v) = X$ for any non-negative number $\varepsilon \geq 0$. 
\end{lemma}

\begin{proof}
Fix any point $x \in X$, any non-negative number $\varepsilon \geq 0$, and any $T>0$. 
Choose any $\varepsilon' \in [0,\varepsilon]$, any $y \in [x]^1_{(\varepsilon',T)}$, and any $\varepsilon_1>0$. 
Put $\varepsilon'' := \varepsilon' + \varepsilon_1>0$. 
It suffices to show  that $y \sim^1_{(\varepsilon'',T)} x$.
By $y \in [x]^1_{(\varepsilon',T)}$, there are a positive number $r \geq T$ and an $(\varepsilon',T)$-link $(v^t(z))_{t \in [0,r]}$ from $x$ to $y$. 
Then $d(x,z) <  \varepsilon'$ and $d(v^r(z),y) < \varepsilon'$. 
From the continuity of $v^r$, there is a positive number $\delta>0$ such that $v^r(B_\delta(z)) \subseteq B_{\varepsilon_1}(v^r(z))$. 
Put $\varepsilon_2 := \min \{ \delta, \varepsilon_1 \}$. 
Since $z$ is non-wandering, there are a point $z' \in B_{\varepsilon_2}(z)$ and a positive number $s > r$ such that $d(v^s(z'), z) < \varepsilon_2$.  
Then $d(y,v^r(z')) \leq d(y,v^r(z)) + d(v^r(z),v^r(z')) < \varepsilon' + \varepsilon_1 = \varepsilon''$ and $d(v^s(z'), x) \leq d(v^s(z'), z) + d(z,x) < \varepsilon_2 + \varepsilon' \leq \varepsilon''$. 
Therefore, a subset $(v^t(z))_{i=r}^{s}$ is an $\varepsilon''$-link from $y$ to $x$. 
This means that $y \sim^1_{(\varepsilon'',T)} x$. 
\end{proof}

We have the following equivalence. 

\begin{theorem}\label{th:02_flow}
Let $v$ be a flow on a metric space $X$.
The following statements are equivalent: 
\\
{\rm(1)} $\Omega(v) = X$. 
\\
{\rm(2)} $\Omega_{-\varepsilon}(v) = X$ for any number $\varepsilon \geq 0$. 
\\
{\rm(3)} $z \sim^1_{\varepsilon+} x$ for any point $x \in X$, any number $\varepsilon \geq 0$, and any point $z \in [x]^1_\varepsilon$.
\end{theorem}

\begin{proof}
The definition of $\varepsilon$-non-wandering point implies that assertions (2) and (3) are equivalent. 
By Lemma~\ref{lem:equv_minus}, assertions (1) and (2) are equivalent. 
\end{proof}

The previous theorem implies a characterization of wandering domains as follows. 
 
 \begin{corollary}\label{cor:wandering_domain_ch}
The following statements are equivalent for a flow $v$ on a metric space $X$:
\\
{\rm(1)} The flow $v$ has a wandering domain. 
\\
{\rm(2)} There are a point $x \in X$, numbers $T>0$ and $\varepsilon \geq 0$, and a point $z \in [x]^1_\varepsilon$ such that $z \notin \bigcap_{\varepsilon'>\varepsilon} [x]^1_{(\varepsilon,T)}$.
\\
{\rm(3)} There are a point $x \in X$, numbers $T>0$ and $\varepsilon \geq 0$, a point $z \in [x]^1_{(\varepsilon,T)}$, and a number $\varepsilon'>\varepsilon$ such that there is no $(\varepsilon',T)$-link from $x$ to $z$. 
 \end{corollary}

\section{Non-wandering diagrams}

\subsection{Non-wandering diagrams for mappings}

As bifurcation diagrams of dynamical systems, we define the non-wandering diagram to analyze dynamical systems as follows. 

\subsubsection{$\varepsilon$-non-wandering diagram of a mapping}

Let $f \colon X \to X$ be a mapping (resp. continuous mapping) on a metric space $X$ and put $I := (-\infty,-0] \sqcup [0,\infty]$ (resp. $I := \R$). 
As a representation of the non-wandering filtration, we introduce the following definition.

\begin{definition}
The subset $D_{\Omega}(f) := \{(\varepsilon, \Omega_{\varepsilon}(f)) \mid \varepsilon \in I \}$ is called the {\bf non-wandering diagram} of $f$. 
\end{definition}

\subsubsection{Examples of non-wandering diagrams for mappings}

Consider an expanding homeomorphism $f \colon \R \to \R$ by $f_2(x):= 2x$. 
Then $\Omega_{-\varepsilon}(f_2) = \emptyset \subsetneq \{ 0 \} = \Omega_{-0} = \Omega_{0}  \subsetneq [-3 \varepsilon, 3 \varepsilon] = \Omega_{\varepsilon}(f_2)$ for any $\varepsilon > 0$.  
This implies the following singular limit. 

\begin{theorem}\label{th:counter_ex04}
For any $\varepsilon>0$, there is a homeomorphism $f$ on a circle with $\bigcup_{\varepsilon>0} \Omega_{-\varepsilon}(f) \subsetneq \Omega_{-0}(f) = \Omega_{0}(f)$. 
\end{theorem}

A homeomorphism $f_{\mathrm{rep}} \colon \R \to \R$ defined by 
\[
f_{\mathrm{rep}}(x) = 
\begin{cases}
x &  (x \leq 0) \\
2x &  (x > 0) 
\end{cases}
\]
satisfies that $\Omega_{-\varepsilon}(f) = \R_{\leq -\varepsilon } \subseteq \R_{\leq 0 } = \Omega_{-0} = \Omega_{0}  \subseteq \R_{\leq 3 \varepsilon} = \Omega_{\varepsilon}(f_{1/2})$ for any $\varepsilon > 0$.  

Moreover, consider a contraction $f_{1/2} \colon \R \to \R$ by $f_{1/2}(x)= x/2$. 
Then the non-wandering diagram of $f_{1/2}$ satisfies the following: 
\[
D_{\Omega}(f_{1/2}) = (\R_{<0} \times \{ 0 \}) \sqcup \{ (\varepsilon, x) \mid \varepsilon \geq 0, x \in [-3\varepsilon, 3\varepsilon]\}
\]
\begin{figure}[t]
\begin{center}
\includegraphics[scale=0.345]{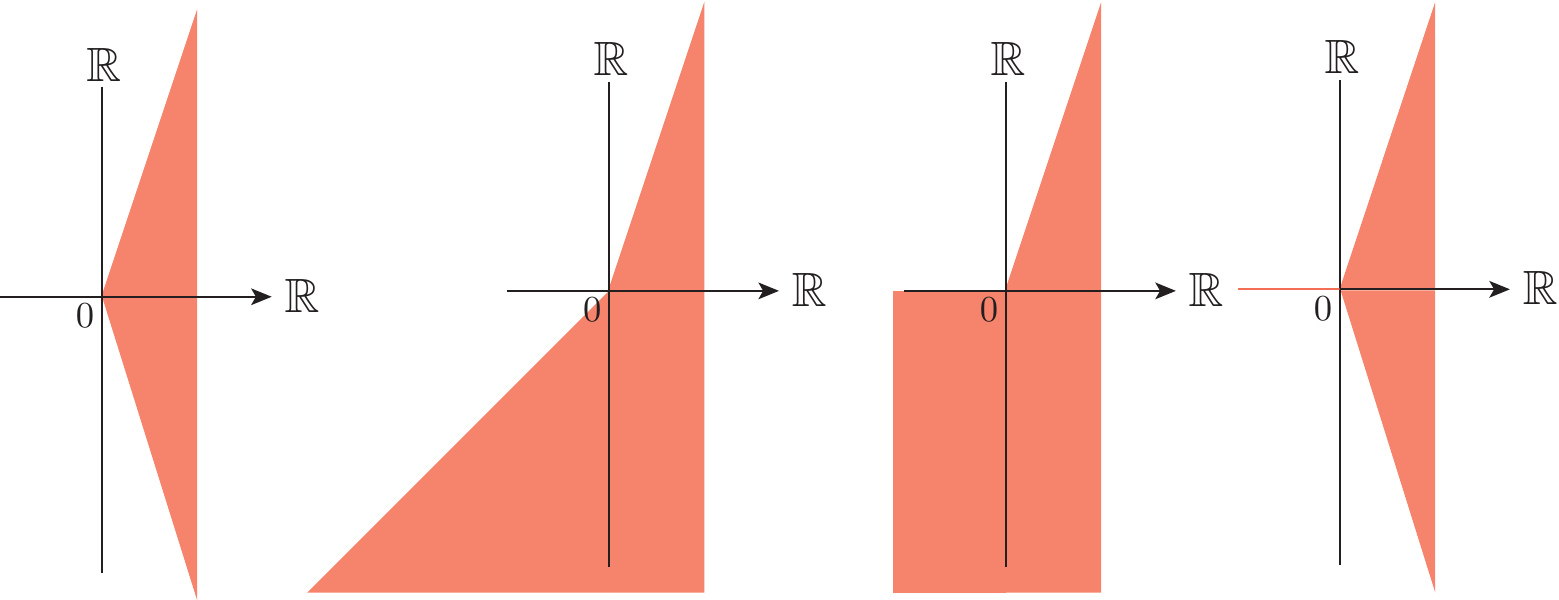}
\end{center} 
\caption{Diagrams $D_{\Omega}(f_2)$, $D_{\Omega}(f_{\mathrm{rep}})$, $D_{\Omega}(f_{\mathrm{att}})$, and $D_{\Omega}(f_{1/2})$.}
\label{fig:ex_diagram}
\end{figure} 

Furthermore, consider the homeomorphism $f_{\mathrm{rep}}$ as above, and a homeomorphism $f_{\mathrm{att}} \colon \R \to \R$ defined as follows: 
\[
f_{\mathrm{att}}(x) = 
\begin{cases}
x &  (x \leq 0) \\
x/2 &  (x > 0) 
\end{cases}
\]
Then the non-wandering diagrams of $f_{\mathrm{rep}}$ and $f_{\mathrm{att}}$ are 
\[
D_{\Omega}(f_{\mathrm{rep}}) = \{ (\varepsilon,x) \mid x \leq \varepsilon \leq 0 \} \sqcup \{ (\varepsilon,x) \mid x \leq 3 \varepsilon, 0 < \varepsilon \}
\]
\[
D_{\Omega}(f_{\mathrm{att}}) = (\R \times \R_{\leq 0}) \sqcup \{ (\varepsilon,x) \mid 0 < x/3 \leq \varepsilon \}
\]
and the non-wandering diagrams of $f_{\mathrm{rep}}$ and $f_{\mathrm{att}}$ are 
\[
D_{\Omega}(f_{\mathrm{rep}}) = \{ (\varepsilon,x) \mid x \leq \varepsilon \}
\]
\[
D_{\Omega}(f_{\mathrm{att}}) = D_{\Omega}(f_{\mathrm{att}}) = (\R \times \R_{\leq 0}) \sqcup \{ (\varepsilon,x) \mid 0 < x/2 \leq \varepsilon \}
\]
as shown in Figure~\ref{fig:ex_diagram}.
The diagrams of these homeomorphisms suggest that the negative parts of the filtrations $D_{\Omega}$ illustrate parts of the behaviors near the non-wandering sets.

\subsection{Non-wandering diagram of a {\rm(}semi{\rm)}flow}

As mappings, we define the non-wandering diagram of (semi)flows to analyze them as follows. 
Let $v$ be a (semi)flow on a metric space $X$.

\begin{definition}
The subset $D_{\Omega}(v) := \{(\varepsilon, \Omega_{\varepsilon}(v)) \mid \varepsilon \in \R \}$ is called the {\bf non-wandering diagram} of $v$. 
\end{definition}

\begin{figure}[t]
\begin{center}
\includegraphics[scale=0.375]{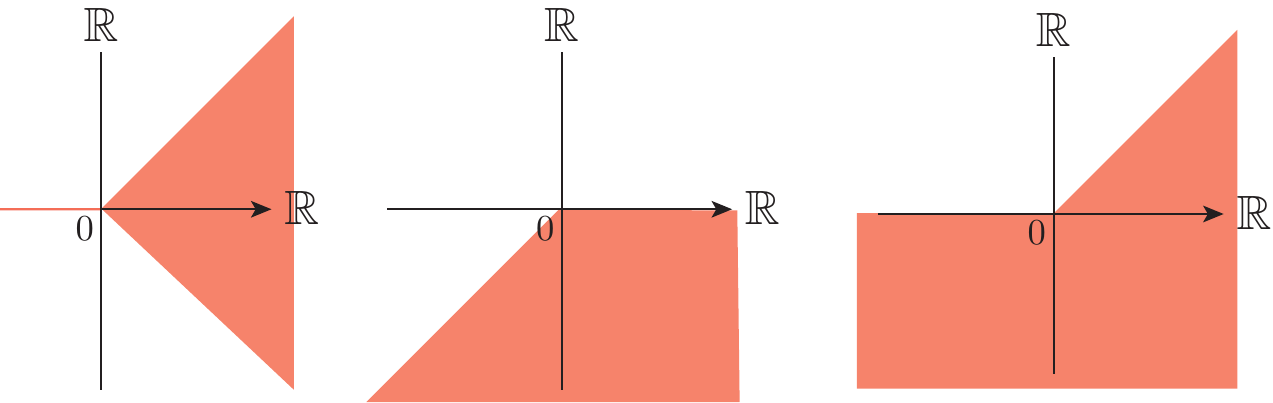}
\end{center} 
\caption{Diagrams $D_{\Omega}(v_Z)$, $D_{\Omega}(v_{\mathrm{rep}})$, and  $D_{\Omega}(v_{\mathrm{att}}) = D_{\Omega}(v_{\mathrm{att}})$.}
\label{fig:csr_diagram_05}
\end{figure}

For instance, consider an attracting flow $v_Z \colon \R \times \R_{\geq 0 } \to \R$ generated by a vector field $Z = -x$. 
Then the non-wandering diagram $D_{\Omega}(v_Z)$ satisfies 
\[
(\R_{<0} \times \{ 0 \}) \sqcup \{ (\varepsilon, x) \mid \varepsilon \geq 0, x \in [-\varepsilon, \varepsilon]\}
\]
as in Figure~\ref{fig:csr_diagram_05}.
Moreover, consider a repelling flow $v_Y \colon \R \times \R_{\geq 0 } \to \R$ generated by a vector field $Y = x$. 
Then the non-wandering diagram $D_{\Omega}(v_Y)$ of $v_Y$ is $\R_{\geq 0} \times \{ 0 \}$. 

In addition, consider flows $v_{\mathrm{rep}}, v_{\mathrm{att}} \colon \R \times \R \to \R$ generated by vector fields $Y_{\mathrm{rep}}, Y_{\mathrm{att}}$, where 
\[
Y_{\mathrm{rep}}(x) = 
\begin{cases}
0 &  (x \leq 0) \\
x &  (x > 0) 
\end{cases}
\]
\[
Y_{\mathrm{att}}(x) = 
\begin{cases}
0 &  (x \leq 0) \\
-x &  (x > 0) 
\end{cases}
\]
are vector fields on $\R$. 
Then the non-wandering diagrams of $v_{\mathrm{rep}}$ and $v_{\mathrm{att}}$ are 
\[
D_{\Omega}(v_{\mathrm{rep}}) = \{ (\varepsilon,x) \mid x \leq \varepsilon < 0 \} \sqcup (\R_{\geq 0} \times \R_{\leq 0}) = \{ (\varepsilon,x) \mid x \leq \varepsilon < 0 \} \sqcup D_{\Omega}(v_{\mathrm{rep}})
\]
\[
D_{\Omega}(v_{\mathrm{att}}) = \{ (\varepsilon,x) \mid 0< x \leq \varepsilon \} \sqcup (\R \times \R_{\leq 0})
\] 
as in Figure~\ref{fig:csr_diagram_05}.
The diagrams of these flows suggest that the negative parts of the filtrations $D_{\Omega}$ correspond to the behaviors near the non-wandering sets.

\section{An example}

The following example demonstrates that the condition $x \sim^1_{-0+} x$ need not imply $x \in \Omega_0(f)$.

\begin{lemma}\label{lem:couter_example01}
There is a continuous mapping $f \colon X \to X$ on a metric space $X$ whose cost function is the metric, and there is a point $p \in X - \Omega(f)$ such that $p \sim^1_{-0 +} p$.
\end{lemma}

\begin{proof}
Put $A : = \{ (1/n,0) \mid n \in \Z_{>0} \}$, $B : = \{ (1/n,1/m) \mid n \in \Z_{>1}, m \in \Z_{>0} \}$ and $X := A \sqcup B \subset \R^2$. 
Define a continuous mapping $f \colon X \to X$ by 
\[
f(x,y) := 
\begin{cases}
\left( \dfrac{1}{n+1},0 \right) &  \left( (x,y) = \left(\dfrac{1}{n},0 \right) \in A \right) \\
(1,0) &  (y = 1) \\
\left( \dfrac{1}{n+1}, \dfrac{1}{m-1} \right)  & \left( y \neq 1 \text{ and } (x,y) = \left(\dfrac{1}{n}, \dfrac{1}{m} \right) \in B \right)
\end{cases}
\]
with respect to the metric induced by the Euclidean metric on $\R^2$. 
By construction, the point $p := (1,0) \in A$ is wandering. 
Since the cost function is the metric $d$, we obtain $O^+(p) = [p]^1_0$. 
For any $z = (1/n,0)  \in O^+(p)$ and any $\varepsilon > 0$, there is a natural number $M \in \Z_{>0}$ with $1/M < \varepsilon$ such that $(1/n, 1/M) \in B_\varepsilon(z)$ and so that $p \in O^+(B_\varepsilon(z)) \subseteq [z]^1_{\varepsilon}$. 
Then $p \in \bigcap_{ z \in O^+(p)}  \bigcap_{\varepsilon > 0} O^+(B_\varepsilon(z)) \subseteq \bigcap_{ z \in O^+(p)}  \bigcap_{\varepsilon > 0} [z]^1_{\varepsilon}$. 
Therefore, the following relations: 
\[
\begin{split}
& p \in \bigcap_{ z \in O^+(p)}  \bigcap_{\varepsilon > 0} [z]^1_{\varepsilon} \text{ with respect to } d 
\\
\Longleftrightarrow \,\, & p \in \bigcap_{ z \in O^+(p)} [z]^1_{0+} \text{ with respect to } d 
\\
\Longleftrightarrow \,\, & p \in [z]^1_{0+} \text{ for any } z \in O^+(p) \text{ with respect to } d 
\\
\Longleftrightarrow \,\, & z \sim^1_{0 +} p \text{ for any } z \in O^+(p) = [p]^1_0 \text{ with respect to } d 
\\
\Longleftrightarrow \,\, & z \sim^1_{0 +} p \text{ for any } z \in [p]^1_0 \text{ with respect to } d 
\\
\Longleftrightarrow \,\, & p \sim^1_{-0+} p \text{ with respect to } d 
\end{split}
\]
This completes the assertion. 
\end{proof}

\bibliographystyle{abbrv}
\bibliography{yt20211011}

\end{document}